\documentclass{article}
\usepackage[T2A]{fontenc}
\usepackage[cp1251]{inputenc}
\usepackage[english,russian]{babel}
\usepackage[tbtags]{amsmath}
\usepackage{amsfonts,amssymb,mathrsfs,amscd}
\usepackage{verbatim}
\usepackage{eucal}
\usepackage{graphicx}
\usepackage{enumerate}
\usepackage[pdftex,unicode,colorlinks,linkcolor=blue,
citecolor=red,bookmarksopen,pdfhighlight=/N]{hyperref}
\usepackage[hyper]{faa-1-arxiv}

\JournalName{%%Функциональный анализ и его приложения
}
\numberwithin{equation}{section}
\theoremstyle{plain}
\newtheorem{maintheo}{\bf Основная теорема}
\newtheorem{theoglue}{\bf Теорема о склейке}
\newtheorem{theorem}{Теорема}
\newtheorem{theoA}{\bf Теорема A}
\newtheorem{lemma}{\bf Лемма}[section]
\newtheorem{propos}{\bf Предложение}[section]

\newtheorem{PLf}{\bf Формула Пуанкаре\,--\,Лелона}
\theoremstyle{definition}
\newtheorem{definition}{Определение}

\newtheorem{proof}{Доказательство}

\newtheorem{example}{Пример}[section]
\newcommand{\const}{{\rm const}}

\renewcommand{\leq}{\leqslant}
\renewcommand{\geq}{\geqslant}
\renewcommand{\d}{\,{\rm d}}
\def\RR{\mathbb R}
\def\CC{\mathbb C}
\def\NN{\mathbb N}

\DeclareMathOperator{\clos}{clos}
\DeclareMathOperator{\Int}{int}
\DeclareMathOperator{\Har}{har}
\DeclareMathOperator{\Hol}{Hol}

\DeclareMathOperator{\Zero}{Zero}
\DeclareMathOperator{\sbh}{sbh}
\DeclareMathOperator{\dom}{dom}
\DeclareMathOperator{\dsbh}{\text{$\delta${\rm -sbh}}}
\DeclareMathOperator{\supp}{supp}
\DeclareMathOperator{\comp}{c}
\DeclareMathOperator{\Conc}{Conc}
\DeclareMathOperator{\gm}{gm}
\DeclareMathOperator{\reg}{reg}
\begin{document}

\title{К распределению нулевых множеств голоморфных функций}

\author[B.\,N.~Khabibullin]{Б.\,Н.~Хабибуллин}
\address{Башкирский государственный университет}
\email{khabib-bulat@mail.ru}
%%второй автор
\author[A.\,P.~Rozit]{А.\,П.~Розит}
\address{Башкирский государственный университет}
\email{rozit@mail.ru}

\date{14.11.2016}
\udk{517.55 : 517.574 : 517.987.1}

\maketitle

\begin{fulltext}

\begin{abstract}
Пусть   $M$ --- субгармоническая функция с мерой Рисса $\nu_M$ в области $D$ из  $n$-мерного комплексного евклидова пространства  $\mathbb C^n$, $f$ --- ненулевая голоморфная в $D$ функция,  
$|f|\leq \exp M$ на $D$ и функция $f$ обращается в нуль на множестве  ${\sf Z}\subset D$. 
Тогда   ограничения на рост меры Рисса $\nu_M$ функции $M$ вблизи границы области $D$ влекут за собой определённые  ограничения на  размеры или площадь/объем множества $\sf Z$.
Количественная форма исследования этого явления даётся в субгармоническом обрамлении. 
\end{abstract}
	
\markright{К распределению нулевых множеств голоморфных функций}

\footnotetext[0]{Работа выполнена при финансовой поддержке РФФИ (грант №\,16-01-00024).}

%%\tableofcontents

\section{Введение}\label{s1}

\subsection{О предшествующих результатах}\label{111}
 Многомерные результаты об описании нулевых множеств голоморфных функций с ограничениями на  рост их модуля  вблизи границы $\partial D$ области определения $D$ до середины 1990-х гг.  содержатся  в книге \cite{GK}, в  статьях \cite{DH78}--\cite{Khab01},  обзорах \cite[\bf 6.5]{Hen85}, \cite[\S~6]{Shv85}. Так, в случае ограниченной ненулевой голоморфной функции широко известно, что объем/площадь ее нулевого множества ограничен (часть классической  теорема Г.\,М.~Хенкина \,--\,Х.~Скоды), а в случае ненулевой голоморфной функции конечного порядка роста  ее нулевое множество того  же порядка роста (чаcть  теоремы Ш.\,М.~Даутова \,--\,Х.~Скоды).    Случай  целых функций многих переменных достаточно полно освещён в \cite{Ron77}--\cite{Khsur}.  Некоторые субгармонические многомерные результаты  подобного типа получены относительно недавно 
С.\,Ю.~Фаворовым и Л.\,Д.~Радченко в   \cite{FR13}.
Для голоморфных и субгармонических функций  в областях на комплексной плоскости $\CC$ история рассматриваемой тематики достаточно детально изложена в \cite{KhT16}. 
В настоящей работе рассматривается  лишь <<лёгкая часть>> задачи об описании нулевых (под)множеств голоморфных функций с заданной мажорантой:  необходимые условия в виде ограничений на рост <<площади--объёма>> нулевого множества вблизи границы области определения. Но  даётся она  в самой общей форме: для произвольных областей и очень широкого круга  условий как по отношению ограничений на рост модуля голоморфных  функций, так и разнообразных   интегральных ограничений на их нулевые множества. Предшествующий  начальный характерный  результат подобного типа  --- 

\begin{theoA}[{\rm \cite[основная теорема]{KhT16_b}, \cite[следствие 1]{KhT16_c}}] Пусть $K$ --- компакт с непустой внутренностью в области $D\subset \CC^n$, $M\not\equiv -\infty$ --- субгармоническая функция в $D$ с мерой Рисса   $\nu_M$, $v$ --- положительная ограниченная субгармоническая функция в $D\setminus K$, удовлетворяющая условию $\limsup\limits_{D \ni z'\to z} v(z')=0$  для любой точки $z\in \partial D$ такой, что $\int_{D\setminus K} v\, {\rm d} \nu_M<+\infty$.
Если $f$ --- ненулевая голоморфная функция в $D$ с множеством нулей ${\sf Z}= \{z\in D\colon f(z)=0\}$, для которой $\log \bigl|f(z)\bigr|\leq M(z)$ при всех $z\in D$, то  
интеграл 	$\int_{{\sf Z}\setminus K}v\, {\rm d} \sigma_{2n-2}$ конечен
для $(2n-2)$-мерной поверхностной меры, или $(2n-2)$-меры Хаусдорфа, $\sigma_{2n-2}$.
\end{theoA}

Мы развиваем теорему A как в отношении равномерно ограниченных подклассов тестовых функций $v$ (ср. с результатами Б.~Коренблюма и С.~Сейпа о равномерных классах  Бергмана \cite{HKZh} и с \cite[\S~5]{KuKh09}), так и в направлении значительно более общих мажорант $M$ --- произвольных разностей субгармонических функций, в виде которых представим очень широкий класс фунций, простейшие из которых  сами субгармонические, а также  произвольные дважды дифференцируемые  функции на $D$.

Зачастую без упоминания, наряду с материалом отмеченных выше работ, используются сведения из  \cite{HK}--\cite{Klimek}.  В   доказательствах не применяются результаты пока не вышедшей в печати работы \cite{KhT16}, тесно связанной по тематике и методам с настоящей статьёй.  

Авторы глубоко признательны рецензенту за ряд полезных замечаний и поправок.
\subsection{Основные обозначения, определения и соглашения}\label{s1.1}
Этот  п.~\ref{s1.1} можно пропустить и обращаться к нему лишь по мере необходимости.
\subsubsection{Множества, топология, порядок}\label{1_1_1}
 $\NN$ --- множество {\it натуральных чисел,\/} $\NN_0:=\{0\}\cup \NN$, а $\RR=\RR^1$ и  $\CC=\CC^1$
 --- множества соответственно всех {\it  вещественных\/} и  {\it комплексных чисел;\/} 
 $\RR_{+\infty}:=\RR\cup \{+\infty\}$, $\RR_{\pm\infty}:=\{-\infty\}\cup \RR_{+\infty}$,
$\RR^+:= \{x\in \RR\colon x\geq 0\}$, $\RR_{+\infty}^+:= \RR^+\cup \{+\infty\}$.
На этих множествах    отношение порядка  $\leq$ на $\RR$ дополняем неравенствами
 $-\infty \leq x\leq +\infty$ для всех  $x\in \RR_{\pm\infty}$. 
   Для $n,m\in \NN$ аффинные пространства  $\CC^n$ и  $\RR^m$ соответственно над $\CC$ и $\RR$ наделяются стандартной евклидовой нормой-модулем $|\cdot|$.
Полагаем $\RR_{\infty}^m:=(\RR^m)_{\infty}$, $\CC_{\infty}^n:=(\CC^n)_{\infty}$, $\CC_{\infty}:=(\CC^1)_{\infty}$ --- одноточечные  компактификации Александрова;   $|\infty|:=+\infty$. При необходимости $\CC^n$ и 
 $\CC_{\infty}^n$ отождествляем соответственно с $\RR^{2n}$ и  $\RR_{\infty}^{2n}$ (над $\RR$). Далее, когда  возможно, обозначения вводятся и определения даются только для $\RR^{m}$ и   $\RR_{\infty}^{m}$.  Для подмножества  $S\subset \RR_{\infty}^m$ через $\clos S$, $\Int S$, $\partial S$ и  $\complement S:= \RR^m_{\infty}\setminus S$  обозначаем соответственно {\it замыкание,\/} {\it внутренность\/}, {\it границу\/} $S$ в  $\RR^m_{\infty}$  и {\it дополнение\/} $S$ до  $\RR^m_{\infty}$.  (Под)область в $\RR_{\infty}^m$ --- открытое связное  подмножество в $\RR_{\infty}^m$. 
 Для $S_0 \subset S\subset \RR^m_{\infty}$ пишем $S_0\Subset S$, если $\clos S_0$ --- {\it компактное подмножество\/} в $S$ в топологии, индуцированной с $\RR_{\infty}^m$ на $S$. Для   
 $r\in \RR_{+\infty}^+$ и $x\in \RR^m$ полагаем $B(x,r):=\{x' \in \RR^m \colon |x'-x|<r\}$ --- открытый шар радиуса $r$ с центром  $x$, $B(r):=B(0,r)$ и $B(x,+\infty)= \RR^m$;  $B(\infty,r):= \{x\in \RR^m \colon |x|>1/r\}$ и $B(\infty, +\infty):=\RR^m_{\infty}\setminus \{0\}$. Одним и тем же символом $0$ обозначаем, по контексту, число нуль, начало отсчета, нулевой вектор, нулевую функцию, нулевую меру и т.\,п. 
 Для подмножества  $X$ упорядоченного векторного пространства  чисел, функций, мер и т.\,п. с отношением порядка \;$\geq$\; полагаем   	$X^+:=\{x\in X\colon x \geq 0\}$  --- все положительные элементы из $X$; $x^+:=\max \{0,x\}$.  {\it Положительность\/} всюду понимается, в соответствии с контекстом, как $\geq 0$.  

Класс борелевских подмножеств множества $S\subset \RR_{\infty}^m$ обозначаем через $\mathcal B(S)$.

\subsubsection{Функции}\label{1_1_2}  Для произвольной функции $f\colon X\to Y$ 
 допускаем, что не для всех  $x\in X$  определено значение $f(x)\in Y$. Сужение функции $f$ на $X_0\subset X$ обозначаем через  $f\bigm|_{X_0}$.  {\it Функция $f$ расширенная числовая,\/} если ее образ --- подмножество в $\RR_{\pm \infty}$. Для расширенных числовых функций полагаем $\dom f:=f^{-1}(\RR)$.
На расширенные числовые функции  {\it отношение  порядка\/} индуцируется с $\RR_{\pm\infty}$ как  {\it поточечное.} При этом для расширенных числовых функций $f\colon X\to Y$  и $g\colon X_1\to Y$  при $X_0\subset X\cap X_1$  пишем <<{\it $f\leq g$ на\/} $X_0$>> и <<{\it $f=g$ на\/} $X_0$>>, если соответственно  $f\bigm|_{X_0}\leq g\bigm|_{X_0}$ и $f\bigm|_{X_0}= g\bigm|_{X_0}$. 

Для  открытого подмножества   $\mathcal O \subset \RR^m_{\infty}$ через  $\Har ({\mathcal O})$ обозначаем  векторное  пространство над $\RR$ {\it гармонических\/} (аффинных  при $m=1$) в ${\mathcal O}$ функций; $\sbh ({\mathcal O})$ --- выпуклый конус над $\RR^+$ {\it субгармонических\/} (выпуклых при $m=1$) в ${\mathcal O}$ функций. 		 Функцию, {\it тождественно равную\/} $-\infty$ на ${\mathcal O}$, обозначаем символом $\boldsymbol{-\infty}\in \sbh ({\mathcal O})$; $\sbh_*(\mathcal O ):=\sbh (\mathcal O)\setminus \{\boldsymbol{-\infty}\}$. 	Для    $\mathcal O \subset \CC^n_{\infty}$ через  $\Hol ({\mathcal O})$  обозначаем  векторное  пространство над $\CC$ {\it голоморфных\/}  функций на ${\mathcal O}$.   При\footnote{При $m=1$ случай $\infty\in \mathcal O\subset \RR_{\infty}$ в настоящей статье  не обсуждается.} $m\geq 2$ функция $u$ голоморфная или (суб)гармоническая в открытой  окрестности  $\mathcal O$ точки $\infty$, если при {\it преобразовании инверсии с центром в нуле\/} (и в точке $\infty$)
 \begin{equation}\label{stK}
\star \colon x\mapsto x^{\star}:= \begin{cases}
0\quad&\text{при $x=\infty$},\\
\frac{1}{|x|^2}\,x\quad&\text{при $x\neq 0,\infty$},\\
\infty\quad&\text{при $x=0$}
\end{cases}
 \end{equation}
преобразование Кельвина $u^{\star}$ функции $u$, определяемое равенствами\footnote{Метка-ссылка над знаками (не)равенства, включения, или, более общ\'о, бинарного отношения и т.\,п., означает, что данное соотношение как-то связано с  отмеченной ссылкой.}  
$u^{\star}(x^{\star})\overset{\eqref{stK}}{=}|x|^{m-2} u(x)$, $x\in \mathcal O\subset \RR_{\infty}^m$,
такая же, но уже в окрестности $\mathcal O^{\star}:=\{x^{\star}\colon x\in \mathcal O\}$ точки $0$. Инверсия $\star_{x_0}$ с центром $x_0\in \RR^m$ определяется как суперпозиция сдвига $x\mapsto x-x_0$ с последующей инверсией $\star_0:=\star$. Для собственного  подмножества $S\subset \RR_{\infty}^m$ класс $\sbh (S)$ состоит из сужений на $S$ функций, субгармонических в некотором, вообще говоря, своём открытом множестве $\mathcal O\subset \RR_{\infty}^m$, включающем в себя $S$; \; $\sbh^+(S):= \bigl(\sbh(S)\bigr)^+$.   Аналогично для $\Har (S)$. Для ${u}\in \sbh ( S)$ полагаем ${(-\infty)}_{u}(S):=\{z\in S \colon {u}(z)=-\infty\}$ --- это {\it $(-\infty)$-мно\-ж\-е\-с\-т\-во\/} функции ${u}$  в $S$ \cite[3.5]{Rans}, где зачастую пишем просто ${(-\infty)}_{u}$, не указывая $S$. 

Через $\const_{a_1, a_2, \dots}\in \RR$ обозначаем постоянные, вообще говоря, зависящие от $a_1, a_2, \dots$ и, если не оговорено противное, только от них;  $\const_{\dots}^+\geq 0$.
\subsubsection{{Меры}}\label{1_1_3}  Далее $\mathcal M (S)$ --- класс {\it борелевских вещественных мер\/} на подмножествах из $S\in \mathcal B ( \RR^m_{\infty})$ со значениями в $\RR_{\pm\infty}$,       иначе называемых  {\it зарядами;\/}   $ \mathcal M_{\comp}(S)$  --- подкласс мер в $\mathcal M (S)$ с компактным {\it носителем\/} $\supp \nu\Subset S$, $\mathcal M^+ (S):=\bigl(\mathcal M (S)\bigr)^+$. Для заряда $\mu \in \mathcal M (S)$ через $\mu^+$, $\mu^-:=(-\mu)^+$ и $|\mu|:=\mu^++\mu^-$ обозначаем соответственно его {\it верхнюю, нижнюю и полную вариации.\/} Мера $\mu \in  \mathcal M (S)$ {\it сосредоточена\/} на измеримом по мере $\mu$ подмножестве $S_0\subset S$, если $\mu (S')=\mu (S'\cap S_0)$ для любого измеримого по мере  $\mu$ подмножества $S'\subset S$. Совокупность всех подмножеств из $S$, на которых сосредоточена мера 
$\mu \in \mathcal M(S)$, обозначаем через $\Conc \mu$.  Очевидно, $\supp \mu \in \Conc \mu$.
Для $x\in \RR_{\infty}^m$ и $0<r\in \RR_+$ полагаем $\mu(x,r):=\mu\bigl(B(x,r)\bigr)$.
Меру Рисса функции ${u}\in \sbh ({\mathcal O})$, $\mathcal O\subset \RR_{\infty}^m$,  чаще всего будем обозначать как 
\begin{equation}\label{df:cm}
\nu_u:= \frac{1}{s_{m-1}} \,\Delta u\in \mathcal M^+(\mathcal O),  \quad\text{где } s_{m-1}:=\frac{2\pi^{m/2}\max\{1, m-2\}}{\Gamma(m/2)}
\end{equation}
--- площадь $(m-1)$-мерной единичной сферы $\partial B(1)$  в $\RR^m$, домноженная на $\max\{1, m-2\}$, $\Delta$ --- {\it оператор Лапласа,\/}  действующий в смысле теории распределений, или обобщённых функций, а $\Gamma$ --- гамма-функция. Такие меры $\nu_u$ --- меры Радона, т.\,е. определяют линейный положительный непрерывный и ограниченный  функционал на пространстве $C_0(\mathcal O)$ непрерывных финитных функций на $\mathcal O$. В частности, $\nu_u(S)<+\infty$ для каждого измеримого  по $\nu_u$ подмножества $S\Subset \mathcal O$. При $u=\boldsymbol{-\infty}\in \sbh(\mathcal O)$ по определению   $\nu_{\boldsymbol{-\infty}}(S):=+\infty$ для $S\subset \mathcal O$.
 
  Через $\lambda_m\in \mathcal{M}^+ (S)$ обозначаем сужения {\it меры Лебега\/} на собственные  борелевские подмножества  $S\subset \RR^m_{\infty}$, где  при $\infty\in S$ предварительно используем инверсию $\star_{x_0}$ с центром $x_0\notin S$; $\delta_x\in \mathcal{M}^+ (S)$ --- {\it мера Дирака в точке\/} $x\in S$, т.\,е. $\supp \delta_x=\{x\}$ и $\delta_x\bigl(\{x\}\bigr)=1$. В обозначении меры Лебега индекс $m$ часто будем опускать. Для $p\in \RR^+$ через $\sigma_p$ обозначаем {\it $p$-мерную (внешнюю) меру Хаусдорфа,\/} или {\it $p$-меру Хаусдорфа,\/} в собственных  открытых подмножествах $\mathcal O\subset \RR_{\infty}^m$ (при $\infty\in \mathcal O$ снова используем инверсию с центром $x_0\notin \mathcal O$). 
 В настоящей работе $p$-мера Хаусдорфа используется лишь при  целом $p{\in} \NN_0$:
  \begin{equation}\label{df:spb}
    \sigma_p(S):=b_p \lim_{r\to 0} \inf \biggl\{\sum_{j\in \NN}r_j^p\,\colon\,  
  S\subset \bigcup_{j\in \NN}B(x_j,r_j), \; 0\leq r_j<r\biggr\}, 
    \quad  b_0:=1, \quad b_p\overset{\eqref{df:cm}}{:=}\frac{s_{p-1}}{p\max\{1, p-2\}}
    \end{equation}
при $p\in \NN$  --- нормирующий множитель, равный объему единичного шара $B(1)$ в $\RR^p$.  В такой нормировке при $p=0$ для любого подмножества $S\subset  \RR_{\infty}^m$   его $0$-мера Хаусдорфа $\sigma_0(S)$ равна мощности $S$, т.\,е. числу точек в $S$, а при $p=m$ имеем   $\sigma_m=\lambda_m$, 
$s_{m-1}\overset{\eqref{df:cm}}{=}\max\{1,m-2\}\,\sigma_{m-1}\bigl(\partial B(1)\bigr)$, $b_p\overset{\eqref{df:spb}}{=}
\sigma_p\bigl(B(1)\bigr)$. 
   В $\CC_{\infty}^n$, отождествленном с $\RR_{\infty}^{2n}$, в настоящей работе будем использовать только $(2n-2)$-меру Хаусдорфа $\sigma_{2n-2}$, которую часто называют $(2n-2)$-мерной поверхностной мерой. При этом    в \eqref{df:cm} и в \eqref{df:spb}
   \begin{equation}\label{df:sbm}
s_{2n-1}\overset{\eqref{df:cm}}{=} \frac{2\pi^n \max\{1,2n-2\}}{(n-1)!}  \, , \quad b_{2n-2}
\overset{\eqref{df:spb}}{:=}\frac{\pi^{n-1}}{(n-1)!}\, .   
   \end{equation} 

\subsubsection{{Нули  голоморфных функций\/} {\rm \cite[\S~11]{Ron77}--\cite{Khsur}, \cite[гл.~1]{Chi}}}\label{zf}
Пусть $D$ ---   подобласть в $\CC_{\infty}^n$, $0\neq f\in \Hol(D)$. {\it Дивизором нулей} функции $f$ называем функцию  $\Zero_f\colon D \to \NN_0$, равную кратности нуля функции $f$ в каждой точке $z\in D$. 
  Для $f=0\in \Hol (D)$ по определению $\Zero_0\equiv +\infty$ на  $D$. Далее в этом пп.~\ref{zf} всюду $0\neq f\in \Hol (D)$. Дивизор нулей $\Zero_f$ полунепрерывен сверху в $D$.  Носитель $\supp \Zero_f$ --- главное аналитическое множество  чистой размерности $n-1$ над $\CC$  и размерности $2n-2$ над $\RR$, для которого $\reg \supp \Zero_f$ --- множество регулярных точек, и всегда  $\sigma_{2n-2} (\supp \Zero_f\setminus \reg \supp \Zero_f)=0$. Пусть   $\reg \supp \Zero_f=\cup_j {\sf Z}_j$ --- представление в виде объединения не более чем счётного числа связных компонент ${\sf Z}_j$, $j=1,2, \dots$. Тогда семейство $\{{\sf Z}_j\}$   локально конечно в $D$, т.\,е. каждое подмножество $S\Subset D$ пересекается лишь с конечным числом компонент ${\sf Z}_j$. Дивизор нулей  $\Zero_f$ постоянен на каждой компоненте ${\sf Z}_j$, т.\,е. однозначно определено значение $\Zero_f({\sf Z}_j)$ для каждого $j\in \NN$.    
  Каждому дивизору нулей $\Zero_f$  сопоставляем  {\it считающую меру нулей\/}   $n_{\Zero_f}\in \mathcal M^+ (D)$, 
определяемую  как мера Радона равенствами $n_{\Zero_f}(\varphi ):=:\int \varphi  \d n_{\Zero_f} \overset{\eqref{df:spb}}{:=}\int \varphi \Zero_f \d \sigma_{2n-2}$ по всем $\varphi \in C_0(D)$, или эквивалентно,  как  борелевская мера на $D$ по правилу $n_{\Zero_f}(B)\overset{\eqref{df:spb}}{=}\sum\limits_{j} \Zero_f({\sf Z}_j)  \sigma_{2n-2}(B\cap {\sf Z}_j)$ для всех $B\in \mathcal B ( D)$.

\begin{PLf}[\cite{Lel}]
Пусть $D\neq \varnothing$ --- собственная подобласть в $\CC_{\infty}^n$, $n\in \NN$, $0\neq f\in \Hol (D)$. Для   меры Рисса $\nu_{\log |f|}$ функции $\log |f|\in \sbh_*(D)$
имеем равенства
\begin{equation}\label{nufZ}
\nu_{\log |f|}\overset{\eqref{df:cm}}{=} \frac{1}{s_{2n-1}}\Delta \log |f| \overset{\eqref{df:sbm}}{=}
\frac{(n-1)!}{2\pi^n \max\{1,2n-2\}}\Delta \log |f| {=}n_{\Zero_f}.
\end{equation}
\end{PLf}
Функция  ${\rm Z}\colon D\to \RR^+$  ---  {\it поддивизор нулей  для\/} $f\in \Hol (D)$, если  ${\rm Z}\leq \Zero_f$ на $D$. Очевидно,  для $f\in  \Hol (D)$ ее дивизор нулей  --- поддивизор нулей для $f$. 
Интегралы по положительной мере с подынтегральными функциями, содержащими поддивизор,  всюду, вообще говоря, понимаем как верхние интегралы.
\subsubsection{{Область $D$ и ее компактификация}}\label{ns} 
В настоящей статье 	далее всюду 
\begin{equation*}
\boxed{D\neq \varnothing \text{ {\sffamily{ --- собственная подобласть в $\RR_{\infty}^m\neq D$ или в $\CC_{\infty}^n\neq D$,}}}}
\end{equation*}
а при $m=1$ еще  и  $\infty \notin \clos D$, т.\,е. $D\subset \RR$ --- открытый	 интервал конечной длины.
Будет  использоваться одноточечная компактификация Александрова такой области  $D$, которую будем обозначать через $D_{\partial}:=D\cup \{\partial D\}$, где добавляемую {\it точку\/} $\partial D$ также называем границей области $D$. В частности, выше уже возникали 
$(\RR^m)_{\partial}=\RR^m_{\infty}$ и $(\CC^n)_{\partial}=\CC^n_{\infty}$.
  Базой открытых  (замкнутых)  окрестностей  этой точки $\partial D\in D_{\partial}$ служат всевозможные множества $(\partial D)\cup (D\setminus S)$, где $S$ пробегает множество всех компактов $S$
из $D$ (соответственно открытых подмножеств $S\Subset D$). 
Далее фраза  <<{\it вблизи $\partial D$}>> означает <<{\it в некоторой проколотой окрестности точки $\partial D\in D_{\partial}$}>>, или <<{\it на множестве вида $D\setminus S$, где $S\Subset D$}>>. 
Для постоянных $\const_{m,n, D, \dots}$ зависимость от размерностей $m,n$ и области $D$ не указываем и не обсуждаем.
%%%\newpage

\section{Следствия из основной теоремы}\label{CorMT}
\subsection{Тестовые субгармонические функции}\label{tsf} 
Для расширенной числовой функции $v$, определенной вблизи\footnote{Здесь можно, а возможно, и нужно использовать понятие множества  ростков функций, определённых вблизи  $\partial D$ изнутри $D$, т.\,е. функций в проколотой окрестности точки $\partial D\in D_{\partial}$, с соответствующими отношением эквивалентности <<{\it совпадают вблизи $\partial D$}>> и факторизацией по этому отношению эквивалентности.} $\partial D$, запись    $\lim\limits_{\partial D}v =y$
\text{($v$ {\it стремится к $y\in Y$ в $\partial D \in D_{\partial}$\/})} будет означать, что $\lim\limits_{D\ni x\to \partial D}v(x)=y$ в $D_\partial$.  В обычной трактовке границы $\partial D$ как замкнутого подмножества {\it компактного пространства\/}  $\RR^m_{\infty}$, это  	эквивалентно равенствам  $\lim\limits_{D\ni x'\to x}v(x') =y$ {\it для всех точек\/} $x\in \partial D$. Для  меры $\mu$ и функции $v$ вблизи $\partial D$ будем использовать обозначение 
\begin{equation}\label{in_p}
\int^{\partial D}v \d \mu \quad \text{для } \quad \int_{D\setminus S}  v \d \mu \quad\text{с каким-либо борелевским подмножеством}\quad S\Subset D,
\end{equation}
когда важна лишь конечность (сходимость) интеграла,  --- аналог распространённой записи несобственных интегралов  $\int^{\infty}$ без указания нижнего предела интегрирования. 

Аналог основных финитных положительных функций теории распределений Л.~Шварца, или обобщённых функций,  через которые определяется отношение порядка на распределениях или мерах/зарядах, в рамках тематики настоящей  работы  предоставляет
\begin{definition}\label{testv} 
	{\it Субгармоническую положительную функцию $v\geq 0$ вблизи   $\partial D\in D_{\partial}$\/} при
$\lim\limits_{\partial D}v{=}0 $ называем {\it тестовой функцией для $D$ вблизи\/ $\partial D\in D_{\partial}$,\/}
	а класс всех таких  функций обозначаем через $\sbh_0^+(\partial D)$.  Для $b\in \RR^+$ и замкнутого в $\RR_{\infty}^m$ подмножества $S\Subset D$, т.\,е. при $S=\clos S\Subset D$,  подклассы тестовых функций 
	\begin{subequations}\label{sbh+}
\begin{align}
\sbh_0^+(D\setminus S;\leq b)&:=\Bigl\{v\in \sbh^+(D\setminus {S})\colon 
	\sup_{D\setminus {S}}v \leq b, \; 	\lim_{\partial D}v
	{=}0 \Bigr\} ,
	\tag{\ref{sbh+}b}\label{sbh+b}
	\\
	\sbh_0^+(D\setminus S;<+\infty)&:=\bigcup_{b\in \RR^+}\sbh_0^+(D\setminus S;\leq b),
	\tag{\ref{sbh+}S}\label{sbh+S}\\
	\intertext{вырождающиеся при $S=\varnothing$ в одноточечное множество $\{0\}$, порождают}
\sbh_0^+(\partial D)&=\bigcup_{S\Subset D}\sbh_0^+(D\setminus S;<+\infty).
\tag{\ref{sbh+}$\partial$}\label{sbh+par}
\end{align}
	\end{subequations}
	\end{definition}

\subsection{Основные результаты для субгармонической мажоранты $M$}\label{ss:sf}
Приведённые в этом п.~\ref{ss:sf} результаты --- следствия   основной теоремы   из \S~\ref{S2} в случае $D\subset \CC_{\infty}^n$, $n\in \NN$.

\begin{theorem}[{\rm (индивидуальная)}]\label{cor:1}  Пусть 	$ M\in \sbh_*(D)$ с мерой Рисса $\nu_M$ и 
для расширенной числовой функции $w$, определённой вблизи $\partial D$, 
$\int^{\partial D} w \d \nu_M \overset{\eqref{in_p}}{<} +\infty$,
 $0\neq f\in \Hol (D)$ и $|f|\leq \exp M$ на $D\subset \CC_{\infty}^n$.
Тогда для  любой функции $v\overset{\eqref{sbh+par}}{\in} \sbh_0^+(\partial D)$ при  ограничении
$v\leq w $ на каком-либо борелевском множестве из $\Conc \nu_M$ вблизи  $\partial D$
в обозначении \eqref{in_p} имеет место соотношение
$ \int^{\partial D} v \, {\rm Z} \d \sigma_{2n-2}\overset{\ref{zf}}{<} +\infty$
для любого поддивизора ${\rm Z}\overset{\text{\rm \ref{zf}}}{\leq} \Zero_f$.
\end{theorem}
Всюду далее в п.~\ref{ss:sf} предполагаем, что $\varnothing \neq \Int S\subset S=\clos S\Subset D$.  
Для  
\begin{equation}\label{c:w}
w\colon  D\setminus S\to \RR_{+\infty}^+, \quad b:=\sup_{x\in \partial  S}\limsup_{D\setminus {S} \ni x'\to x} w(x')<+\infty 
\end{equation}
{\it наибольшую миноранту относительно конуса\/} $\sbh_0(D\setminus S):=
\bigl\{v\in\sbh ( D\setminus  S)\colon\; \lim\limits_{\partial D}v=0\bigr\}$ определим как функцию
$\gm w:=\sup \{v\overset{}{\in}\sbh_0 ( D\setminus  S)\colon\; v\leq w \text{ на } D\setminus  S\}$.
При этом через  $\gm^*w$ обозначаем {\it полунепрерывную сверху регуляризацию\/} функции $\gm w\colon D\setminus S\to \RR_{+\infty}^+$.

\begin{theorem}[{\rm (индивидуальная)}]\label{th:2}	Пусть 	$ M\in {\text{\rm sbh}}(D)$ с мерой Рисса $\nu_M$, $M\not\equiv -\infty$,  $w$ ---  полунепрерывная сверху функция из\/  \eqref{c:w},  удовлетворяющая одному из двух\footnote{Простые и наглядные достаточные условия регулярности области геометрического характера см., например, в \cite[Теорема 2.11]{HK}, а для $m=2$, т.\,е. в $\CC$,  --- в \cite[4.2]{Rans}.} условий: 
  \begin{equation}\label{key}
 {\rm (i)} \; \lim\limits_{\partial D} w = 0 \quad \text{ или \quad  {\rm (ii)}  $D$ --- регулярная область (для задачи Дирихле).}
  \end{equation}
Пусть также $\int_{D\setminus S} w \d \nu_M< +\infty$.
	Если $0\neq f\in \Hol(D)$ и $|f|\leq \exp M$ на $D$, то для любого поддивизора ${\rm Z}\overset{\text{\rm \ref{zf}}}{\leq} \Zero_f$ выполнено соотношение $\int_{{\rm Z} \setminus S} (\gm^* w) \, {\rm Z} \d \sigma_{2n-2}< +\infty$.
\end{theorem}
При исследовании нулевых множеств в <<жестких>> весовых классах голоморфных функций с более или менее явными мажорантами $M$ возникала необходимость и в равномерных оценках сверху интегралов от  функций тестового характера из некоторого класса по нулевому множеству голоморфной функции. Наиболее рельефно такой подход реализован для равномерных  пространств Бергмана и им подобных в единичном круге  \cite[гл.~4]{HKZh},  \cite[\S\S\,5,6]{KuKh09} и в единичном шаре \cite{BM}.  Следующий результат даёт подобные оценки в общей ситуации.  

\begin{theorem}[{\rm (равномерная)}]\label{th:1} 
Пусть     $z_0\in \Int S$,  $M(z_0)\neq -\infty$ для   $M\in \sbh_* (D)$ с зарядом Рисса $\nu_M$,   
 $b\in \RR^+$. Тогда найдутся постоянные $C:=\const_{z_0,S, b}^+>0$
и $\overline C_M:=\const_{z_0,S, M}^+$, с которыми для любой ненулевой функции $f\in\Hol(D)$ с ограничением $|f|\leq \exp M$ на $D$, а также  для любой тестовой функции $v\overset{\eqref{sbh+b}}{\in}  \sbh_0^+(D\setminus S;\leq b) $ выполнено неравенство
\begin{equation}\label{iq:ufu}
\int_{D\setminus {S}}v\, \Zero_f  \d \sigma_{2n-2}\overset{\ref{zf}}{\leq} 
\int_{D\setminus {S}}  v \,d {\nu}_M	-C \log \bigl|f(z_0)\bigr| +C\,\overline C_M.
\end{equation}
\end{theorem}  

\section{Основная теорема}\label{S2}

\subsection{$\delta$-субгармонические функции и  основная теорема}\label{2_3} 
Следуя {\rm \cite{Ar_d}--\cite{Gr}}, 
при  каждом $m\in \NN$  определим функцию $h_1(t)=t$,
 $h_2(t):=\log |t|$ и $h_m(t):=-|t|^{2-m}$ при $m\geqslant 3 $.  Наряду с функцией $\boldsymbol{-\infty}$ рассматриваем и функцию $\boldsymbol{+\infty}$, тождественно равную $+\infty$. Функция $M\colon D\to \RR_{\pm\infty}$ в $D\subset \RR_{\infty}^m$ 
{\it тривиальная $\delta$-суб\-г\-ар\-м\-о\-н\-и\-ч\-е\-с\-кая,\/} если $M=\boldsymbol{-\infty}$ или $M=\boldsymbol{+\infty}$ на $D$, и   
{\it нетривиальная $\delta$-суб\-г\-ар\-м\-о\-н\-и\-ч\-е\-с\-кая функция с  зарядом  Рисса\/} $\nu_M\in \mathcal M(D)$, если выполнены  следующие три условия-соглашения.
\begin{enumerate}
\item\label{da} {\it Существуют   $ u_1, u_2\in \sbh (D)\setminus \{\boldsymbol{-\infty}\}$ с мерами Рисса  $\nu_{u_1}, \nu_{u_2}\in \mathcal M^+(D)$, для которых $M(x):=u_1(x)-u_2(x)\in \RR$ при $x\notin (-\infty)_{u_1}\cup (-\infty)_{u_2}$.\/} Для $\nu_M:=\nu_{u_1}- \nu_{u_2}\in \mathcal M(D) $  однозначно определено разложение Хана--Жордана $\nu_M:=\nu_M^+-\nu_M^-$. 
\item\label{dad} {\it Определяющее множество ${\dom M}\subset D$ --- это множество точек $x\in D$, для каждой из которых при некотором  $r_x>0$ конечен один из интегралов}
\begin{equation}\label{dMa}
\left(\,\int_0^{r_x}\frac{|\nu_M|(x,t)}{t^{m-1}} \d t<+\infty\right)\Longleftrightarrow
\left(\; \int_{B(x,r_x)} h_m\bigl(|x'-x|\bigr)\d |\nu_M|(x')>-\infty\right).
\end{equation}
В случае $m=1$ всегда $\dom M=D$, т.\,е. \eqref{dMa} выполнено для всех $x\in D$.  {\it При $m>1$ доопределяем функцию $M$ на всех  $x\in \dom M$ через усреднения по сферам $\partial B(x,r)$
\begin{equation*}
		M(x)\overset{\eqref{df:cm}}{=}\lim_{r\to 0}\frac{1}{s_{m-1}r^{m-1}} \int_{\partial B(0,r)} 
M(x+x') \d \sigma_{m-1}(x')	\, \in \RR	\quad \text{для $x\in {\dom M}$},
		\end{equation*}}
Такое доопределение всегда согласуется с предварительными  значениями функции  $M$ в предыдущем п.~\ref{da} на   $\RR_{\infty}^m\setminus 
\bigl((-\infty)_{u_1}\cup (-\infty)_{u_2}\bigr)\subset \dom M $. В частности \cite{HK}--\cite{Rans}, для $u\in \sbh_*(D)$ 	её определяющее множество $\dom u=D\setminus (-\infty)_u$.
\item\label{dai} {$M(x)=+\infty$ при\footnote{В \cite[2]{Gr} положено $M(x)=0$ для $x\notin {\dom M}$, но для наших целей предпочтительнее  $=+\infty$, поскольку функция $M$ у нас играет роль	мажоранты.} $x\in D\setminus {\dom M}$.}
\end{enumerate}

Пункт \ref{da}, следуя \cite[теорема 11]{Ar_d},    
можно заменить, избегая  явного упоминания  субгармонических функций, на пункт
\begin{enumerate}
\item[$1'.$]  {\it $M$ --- локально интегрируемая по мере  Лебега $\lambda_m$ на $D$ функция, обладающая 	свойством: для любой подобласти $D'\Subset D$ найдется постоянная $C'\in \RR^+$, с которой для любой финитной дважды непрерывно дифференцируемой  функции  			$f \colon D'\to \RR$ с носителем $\supp f \subset D'$ имеем 			$\Bigl|\int_{D'} M	\Delta f \, d\lambda_m \Bigr|\leq C' \max_{z\in D'} \bigl|f (z)\bigr|$. При этом $\nu_M\overset{\eqref{df:cm}}{:=}\dfrac{1}{s_{m-1}}\, \Delta M$ в смысле теории распределений.}
\end{enumerate}
При таком подходе    п.~\ref{dad} заменяем  соответственно на пункт
\begin{enumerate}
\item[$2'.$]\label{2st}  {\it Для определяющего множества\/} ${\dom M}\subset D$, определяемого, как и выше, условием конечности интегралов вида  \eqref{dMa}, полагаем 
 \begin{equation*}
M(x)\overset{\eqref{df:spb}}{=}\lim_{r\to 0} \frac{1}{b_mr^m} \int_{B(0,r)} M(x+x') \d \lambda_m(x')
\quad \text{для всех $x\in {\dom M}$},
\end{equation*}
где интеграл с множителем --- {\it усреднение по шару\/} $B(x,r)$, а  $D\setminus {\dom M}$ --- множество нулевой ёмкости \cite[гл.~5]{HK}.  В частности, $\lambda_m (D\setminus {\dom M})=0$.
		\end{enumerate}
		Пункт~\ref{dai} сохраняем или считаем функцию $M$ неопределенной в $D\setminus {\dom M}$.

Класс $\delta$-субгармонических в $D$ функций обозначаем через $\dsbh (D)$, а подкласс нетривиальных ---
$\dsbh_* (D):=\dsbh (D)\setminus \{\boldsymbol{\pm\infty}\}$.
\begin{maintheo}\label{ss:23} 
Пусть\/  $ M\in \dsbh_* (D)$ с зарядом Рисса $\nu_M$ и   определяющим множеством\/   ${\dom M}\subset D$, а также   $\varnothing \neq \Int S\subset S=\clos S \Subset D\subset \RR_{\infty}^m\neq  D$.
Тогда   для любых  точки $x_0\in \Int S\cap {\dom M}$ и числа  $b>0$, а также регулярной области $\widetilde{D}\subset \RR_{\infty}^m$ с функцией Грина $g_{\widetilde{D}}(\cdot , x_0)$ с полюсом в точке $x_0$ при  $S\Subset \widetilde{D}\subset D$  и $\RR_{\infty}^m\setminus \clos \widetilde{D}\neq \varnothing$ 
с числом
\begin{equation}\label{cz0C}
  C:=\const_{x_0,S,\widetilde{D},b}^+:= \frac{b} {\inf\limits_{z\in \partial S}  g_{\widetilde{D}}(x, x_0)}>0,
\end{equation}
 для любой функции $u\in \sbh_* (D)$, удовлетворяющей неравенству 
$u\leq M$ на  $D$,  а также  для любой тестовой функции $v\overset{\eqref{sbh+b}}{\in}  \sbh_0^+(D\setminus S;\leq b) $ выполнено неравенство
\begin{equation}\label{mest}
C u(x_0) 	+\int_{D\setminus S}  v \,d {\nu}_u 		\leq	\int_{D\setminus S}  v \,d {\nu}_M	+\int_{\widetilde{D}\setminus S} v \,d {\nu}_M^-   +C\, \overline{C}_M,
\end{equation}
где для $\overline{C}_M:=	\int_{\widetilde{D}\setminus \{x_0\}} g_{\widetilde{D}}(\cdot, x_0)  \,d {\nu}_M  
		+\int_{\widetilde{D}\setminus S} g_{\widetilde{D}}(\cdot, x_0)  \,d {\nu}_M^-  +M^+(x_0)$ 
  возможно значение $+\infty$, но при $\widetilde{D}\Subset D$ --- это  некоторая постоянная 
$\overline{C}_M:=\const_{x_0,S, \widetilde{D},M}^+<+\infty$.
\end{maintheo}
   
\subsection{Меры и потенциалы Йенсена}\label{i:mpJ}\label{mJ:Ex}  

\begin{definition}[(\cite{Khab01}, \cite{Khab99}, \cite{Khsur}, \cite{KhT16}--\cite{KhT16_c}, \cite{CR}, \cite{Khab03})]\label{df:J} Мера $\mu \in \mathcal M^+(\RR_{\infty}^m)$ ---  {\it мера Йенсена внутри области\/ $D$ в точке\/ $x_0\in D$,\/} если $\mu \in \mathcal M^+_{\comp}(D)$  и 
$u(x_0)\leq \int u\,d \mu$ для всех  $u\in \sbh(D)$.  Класс всех таких мер Йенсена обозначаем через $J_{x_0}(D)$. 
\end{definition}
Очевидно, каждая мера  $\mu \in J_{x_0}(D)$ {\it вероятностная,\/} т.\,е. $\mu (D)=1$. Далее без упоминания используется то, что {\it для любого множества нулевой ёмкости\/} $E\subset D$, в частности, при  $E=(-\infty)_u$, $ u\in \sbh_* (D)$, {\it для меры $\mu \in J_{x_0}(D)$ имеем $\mu \bigl(E\setminus \{x_0\}\bigr)=0$ \cite[следствие  1.8]{CR}.}

 \begin{example}\label{ex:1hm} Пусть $m\geq 2$ и  $\widetilde{D}\neq \varnothing$ --- подобласть в $D$ с неполярной границей,  а также $\RR_{\infty}^m\setminus \clos \widetilde{D}\neq \varnothing$. 
{\it Гармоническая мера\/} $\omega_{\widetilde{D}}(x_0, \cdot)$ {\it для\/} (или относительно) $\widetilde{D}$ {\it в точке\/} $x_0\in \widetilde{D}$ 
{\large(}см. \cite[4.3]{Rans}, \cite[3.6, 5.7.4]{HK}{\large)}
при условии $\widetilde{D}\Subset D$   --- пример меры Йенсена из $J_{x_0}(D)$. Вырожденный случай   --- мера Дирака  $\delta_{x_0}$. 

При $m=1$ гармоническая мера для интервала
 $\widetilde D=(a,b)\Subset \RR$ в точке $x_0\in (a,b)$ --- это мера  $\omega_{(a,b)}(x_0,\cdot )=\frac{x_0-a}{b-a}\delta_a+\frac{b-x_0}{b-a}\delta_b$. 
 Для областей-интервалов  вида $\widetilde D=(a,+\infty)\subset D\subset \RR$ и 
$\widetilde D=(-\infty,a)\subset D\subset \RR$ с $a\in \RR$   
 гармонические  меры не существуют ни в одной точке $x_0\in D$.
 \end{example}

\begin{definition}[{\large(}\cite{Khab01}, 	\cite{Khab99}, \cite{Khsur},  \cite{KhT16}--\cite{KhT16_c}, \cite{Khab03}{\large)}]\label{df:PoJ} Функцию $V\in \sbh^+\bigl(\RR_{\infty}^m \setminus \{x_0\}\bigr)$ 
при $x_0\in D$ называем {\it потенциалом Йенсена  внутри  $D$ с полюсом в  $x_0$,} если выполнены два условия:
\begin{enumerate}[{\rm 1)}]
\item\label{V:f2} {\it найдётся  область $D_V\Subset D$, содержащая  точку $x_0\in D_V$, для которой $V (x)\equiv 0$ при всех   $x\in \RR_{\infty}^m \setminus D_V$, т.\,е. $V\bigm|_{\RR_{\infty}^m \setminus D_V}=0$} ( финитность в $D$);  
\item\label{V:f3} задана  {\it полунормировка в точке\/ $x_0$}, а именно: 
\begin{subequations}\label{nvz}
\begin{align}
\limsup_{x\to x_0}V(x)<+\infty&\quad\text{при $m=1$ и  $x_0\in \RR$},
\tag{\ref{nvz}a}\label{nvz1}
\\
\limsup\limits_{x \to x_0}\dfrac{V(x )}{-\log |x-x_0|}\leq  1 &\quad \text{при $m=2$},
 \tag{\ref{nvz}b}\label{nvz2}\\ 
\limsup\limits_{x\to x_0}|x-x_0|^{m-2}\, V(x )\leq  1  \quad   &\quad\text{при $m\geq 3$ и  $x_0\neq \infty$,} 
\tag{\ref{nvz}c}\label{nvz3}
\\
\limsup\limits_{ x \to \infty}{V(x )}{\leq} 1 &\quad\text{при $m\geq 3$ и $x_0= \infty$.}
\tag{\ref{nvz}d}\label{nvz4}
\end{align}
\end{subequations}
\end{enumerate}
Класс всех  таких потенциалов Йенсена обозначаем через  $PJ_{x_0} (D)$. 
\end{definition}
По определениям \ref{testv} и \ref{df:PoJ}  имеет место очевидное 
\begin{propos}\label{PJ=test}
 Каждый потенциал Йенсена из  $PJ_{x_0} (D)$ --- тестовая функция  для $D$ вблизи\/ $\partial D\in D_{\partial}$, а при  $S\Subset D$ и $\Int S\neq \varnothing$  для любой  точки  $x_0\in \Int S$
справедливы включения
$ PJ_{x_0} (D) \overset{\eqref{sbh+b}}{\subset}  \sbh_0^+\bigl(D\setminus S; \leq \sup_{\partial S} V\bigr) 
\overset{\eqref{sbh+S}}{\subset}  \sbh_0^+\bigl(D\setminus S; <+\infty\bigr)$.
 \end{propos}
\begin{example}\label{exp2}  В условиях примера \ref{ex:1hm}
{\it функция Грина  $g_{\widetilde{D}}(\cdot ,x_0)$} {\it для} (или относительно) $\widetilde{D}${\; \it  с полюсом в точке $x_0\in \widetilde{D}$, продолженная на\/} $\RR_{\infty}^m$ по правилу
{\large(}см. \cite[3.7, 5.7]{HK}, \cite[4.4]{Rans}{\large)}
\begin{equation}\label{gD0r}
	g_{\widetilde{D}}(x, x_0):=\begin{cases}
	\limsup\limits_{\widetilde{D}\ni x'\to x} g_{\widetilde{D}}(x', x_0)\quad &\text{при $x\in \partial \widetilde{D}$},\\
	0 \quad &\text{при $x\in \RR^m_{\infty}\setminus \clos \widetilde{D}$}
	\end{cases}
\end{equation}
   --- пример потенциала  Йенсена из $PJ_{x_0}(D)$  при условии $\widetilde{D}\Subset D$. Вырожденный вариант   --- функция, тождественно равная нулю на $\RR_{\infty}^m\setminus \{x_0\}$. 
 \end{example}
	
	Всюду далее во избежание  технически громоздкого  разбора различных случаев мы избегаем рассмотрения потенциалов Йенсена и функций Грина с полюсом в точке $\infty$. 
	
\begin{definition}\label{df:lp}
{\it Логарифмический потенциал рода\/ $0$  меры\/ $\mu \in \mathcal  M^+_{\comp}(\RR^m_{\infty})$  с полюсом в  $x_0\in \RR^m$} определяем 
		для всех $y\in \RR^m_{\infty}\setminus \{x_0\}$ как функцию \cite[определение 3]{Khab03}
\begin{equation}\label{df:Vmu}
	V_{\mu}(y) := 		\int_{D} \Bigl(h_m\bigl(|y-x|\bigr)-h_m\bigl(|y-x_0|\bigr)\Bigr) \d \mu (x),
\end{equation}
где при $y=\infty$ подынтегральное выражение доопределено значением  $0$. 
\end{definition}

Дальше подробно обсуждаются  только размерности $m\geq 2$, поскольку одномерная ситуация $m=1$  специфична и, на наш взгляд, требует отдельного рассмотрения. В то же время и случай $m=1$ вписывается в нашу общую схему, а  основная теорема и ее следствия из п.~\ref{ss:sf} справедливы и для $m=1$. Отметим для $m\geq 2$ основные взаимосвязи между мерами и потенциалами Йенсена. 
Первая --- следующее утверждение о двойственности.
\begin{propos}[{\cite[предложение 1.4, теорема двойственности]{Khab03}}]\label{pr:1} Отображение  
\begin{equation}\label{con:P}
	\mathcal P \colon J_{x_0}(D)\to PJ_{z_0} (D), \quad 
	\mathcal P (\mu)\overset{\eqref{df:Vmu}}{:=} V_{\mu}, \quad \mu \in  J_{x_0}(D) ,
\end{equation}
--- биекция,  $\mathcal P\bigl(t\mu_1+(1-t)\mu_2\bigr)=t\mathcal P (\mu_1)+(1-t)\mathcal P (\mu_2)$ для всех $t\in [0,1]$,  а 
\begin{equation}\label{eq:mu}
	{{\mathcal P}}^{-1}(V){=}\frac1{s_{m-1}} \Delta  V\Bigm|_{D\setminus \{x_0\}}+
	\left(1-\limsup\limits_{x \to x_0}\dfrac{V(x )}{-h_m\bigl(|x-x_0|\bigr)}\right)\,
	{\delta}_{x_0}\, , \quad V\in PJ_{x_0}(D).
\end{equation}
В частности, для регулярной  области $\widetilde{D}\subset D$ при $x_0\in \widetilde{D}$ --- это классическое равенство
$\mathcal P\bigl(\omega_{\widetilde{D}}( x_0, \cdot )\bigr)=g_{\widetilde{D}}(\cdot, x_0 )$,  $x_0\in \widetilde{D}\subset D$.
\end{propos}

Вторая   ---  это развитие {\it классической  формулы Пуассона\,--\,Йенсена\/} \cite[теорема 3.14]{HK}.
\begin{propos}[{\cite[Предложение 1.2]{Khab03}}]\label{pr:2}
Пусть $\mu \in J_{x_0}(D)$. Тогда для  $u\in \sbh (D)$ с мерой  Рисса
$\nu_u$ при $u(x_0)\neq -\infty$ справедлива расширенная формула Пуассона\,--\,Йенсена
\begin{equation}\label{f:PJ}
	u(x_0) +\int_{D\setminus \{x_0\}} V_{\mu} \d {\nu}_u=\int_{D} u \d \mu  .
\end{equation}
В частности, для ограниченной  области $D$ в $\RR^m$ при $\mu =\omega_D(x_0,\cdot)$ и, как следствие, с
 $V_{\mu}\overset{\eqref{df:Vmu}}{=}g_D(\cdot,x_0)$  --- это обобщённая  формула Пуассона\,--\,Йенсена\/ \cite[теорема 5.27]{HK}.
\end{propos}
\subsubsection{{Продолжение тестовых функций на $D\setminus \{x_0\}$}}\label{spf:241}   

\begin{propos}\label{prop:V} При $\varnothing \neq \Int S\subset S=\clos S\Subset D$ для любых 	
 точки $x_0\in \Int S$ и числа  $b\in \RR^+$,  а также регулярной области
	 	$\widetilde{D}\subset \RR_{\infty}^m$ при\/ $S\Subset \widetilde{D}\subset D$  и $\RR_{\infty}^m \setminus \clos \widetilde{D}\neq \varnothing$ с числом   
\begin{equation}\label{cz0}
 \widetilde{c}:=\const_{x_0,S,\widetilde{D}, b}^+=\frac{1}{b} \inf_{x\in \partial S}  g_{\widetilde{D}}(x, x_0)  >0
\end{equation}
  по любой тестовой функции $v\overset{\eqref{sbh+}}{\in}  \sbh_0^+(D\setminus S;\leq b) $ с мерой Рисса $\mu_v$
при помощи продолженной функции Грина  $g_{\widetilde{D}}(\cdot, x_0)$ можно  построить  функцию 
\begin{equation}\label{cr:Vng}
\widetilde{V}:= 
\begin{cases}
g_{\widetilde{D}}(\cdot, x_0)\quad &\text{на $S\setminus \{x_0\}$},\\
\max \bigl\{g_{\widetilde{D}}(\cdot, x_0) ,\; \widetilde{c}\cdot v\bigr\} \quad &
\text{на $\widetilde{D}\setminus S$}
,\\
\widetilde{c}\cdot v \quad &\text{на $D\setminus \widetilde{D}$},\\
0 \quad &\text{на $\RR^m_{\infty}\setminus {D}$},
\end{cases}
\quad \qquad \widetilde{V}\in \sbh^+\bigl(\RR^m_{\infty}\setminus \{x_0\}\bigr),
\end{equation}
обладающую свойствами
\begin{equation}\label{VVV}
\widetilde{V}\bigm|_{\Int S \setminus \{x_0\}} \in \Har \bigl(\Int S\setminus \{x_0\}\bigr),
\quad  \lim_{\partial D} \widetilde{V}{=}0,
\quad \lim_{ x \to x_0}
\frac{\widetilde{V}(x )}{-h_m\bigl(|x-x_0|\bigr)}\overset{\eqref{nvz2}{-}\eqref{nvz4}}{=}  1.
\end{equation}
\end{propos}

\begin{proof}
Функцию $v\overset{\eqref{sbh+}}{\in} \sbh_0^+(D\setminus S;\leq b)$ 
можно рассматривать как продолженную нулем на $\complement D$, используя для её обозначения  ту же букву $v$. Очевидно,  продолженная функция $v$  положительная  субгармоническая всюду на  $\RR_{\infty}^m\setminus S$ и $v=0$ на $\complement D$. Из  ограничения  $\sup\limits_{z\in  D\setminus S}v\leq b$ в определении \eqref{sbh+b} следует  $\sup\limits_{x\in \partial  S}\limsup\limits_{D\setminus {S} \ni x'\to x} v(x')\leq b$.
Тогда
\begin{equation}\label{ev0}
	\limsup_{(\RR^m_{\infty} \setminus S) \ni x'\to x} v(x')\leq  b\quad 
	\text{для всех $x\in \partial (\RR^m_{\infty} \setminus S)$.}
\end{equation}
Для области $\widetilde{D}$ из условия предложения \ref{prop:V} рассмотрим  продолженную как в \eqref{gD0r} функцию Грина $g_{\widetilde{D}}(\cdot, x_0)$.
Из  свойств функции Грина \cite[3.7]{HK}, \cite[теорема 4.4.3]{Rans} она гармонична и строго положительна 
на $\widetilde D\setminus \{x_0\}$, откуда
$a:=\inf\limits_{x\in \partial (\RR^m_{\infty} \setminus S)}  g_{\widetilde{D}}(x, x_0)>0$.	
Для функции 
$v_0{:=}\frac{b}{a}\,  g_{\widetilde{D}}(\cdot , x_0)  \in \sbh^+\bigl(\RR^m_{\infty}\setminus\{x_0\}\bigr)$
по  свойствам функций Грина 
\cite[теорема 4.4.9]{Rans}, \cite[1.5.1]{HK}
\begin{equation}\label{v0o}
 \lim_{D\ni x\to x_0}\frac{v_0(x)}{-h_m\bigl(|x-x_0|\bigr)}
{=} \frac{b}{a}\,, \quad v_0 = 0 \quad\text{на  $ D\setminus \widetilde{D}$}, \quad 
v_0\bigm|_{\widetilde{D}\setminus \{x_0\}}\in \Har \bigl(\,\widetilde{D}\setminus \{x_0\}\bigr).
\end{equation}
Кроме того, ввиду  \eqref{ev0}-- \eqref{v0o} и по построению функции $v_0$ имеем
\begin{equation}\label{vvv}
	v_0\bigm|_{\partial (\RR^m_{\infty} \setminus S)}\overset{\eqref{v0o}}{\geq} b
	\overset{\eqref{ev0}}{\geq} 
		\sup_{x\in \partial (\RR^m_{\infty} \setminus S)}\; \limsup_{(\RR^m_{\infty} \setminus S) \ni x'\to x} v(x'), 
\quad 	v_0\bigm|_{\complement \widetilde{D}}=0\leq v\bigm|_{\complement  \widetilde{D}} .
\end{equation}

\begin{theoglue}[{\rm (см.  \cite[следствие  2.4.4]{Klimek})}]
Пусть  $\mathcal  O$ и $\mathcal  O_0$ --- открытые множества в $\RR^m_{\infty}$ и 
$\mathcal O \subset \mathcal O_0$. Пусть  $v\in \sbh (\mathcal O)$ и $v_0\in \sbh (\mathcal O_0)$.  Если 
\begin{equation}\label{0vs}
\limsup_{\mathcal O\ni  x'\to x} v(x')\leq v_0(x) \quad \text{для всех точек $x\in \mathcal O_0\cap \partial \mathcal O$},
\end{equation}
то функция
\begin{equation}\label{consv}
\widetilde{v}:=\begin{cases}
\max\{v,v_0\} \quad &\text{на\/ $\mathcal O$},\\
v_0		\quad &\text{на\/ $\mathcal O_0\setminus \mathcal O$},
\end{cases}
\end{equation}
--- субгармоническая на $\mathcal O$, т.\,е.  $\widetilde v\in \sbh(\mathcal O_0)$. 	
\end{theoglue}
 Применим теорему о склейке при
$\mathcal O_0:=\RR_{\infty}^m\setminus \{x_0\}$ и $\mathcal O:=\RR_{\infty}^m\setminus S$
к функции $v_0$  и продолженной  на  $\mathcal O$ функции $v$. Ввиду первого соотношения из \eqref{vvv}  эти функции удовлетворяют условию \eqref{0vs}. По построению
\eqref{consv} и  из \eqref{v0o}--\eqref{vvv}  для построенной функции $\widetilde{v}\in \sbh^+\bigl(\RR_{\infty}^m\setminus \{x_0\}\bigr)$ её конструкцию можно описать через соответствующие сужения более детально:
\begin{equation}\label{v000}
	0\leq \widetilde{v}=
	\begin{cases}
		v_0 \quad &\text{на $S\setminus \{z_0\}$},\\ 
	\max\{v,v_0\}\quad &\text{на $\widetilde{D}\setminus S$},\\
	v \quad &\text{на $D\setminus  \widetilde{D}$},\\
	0 \quad &\text{на $\RR^m_{\infty}\setminus  D$},
	\end{cases}
\quad \qquad \widetilde{v}\in \sbh^+\bigl(\RR^m_{\infty}\setminus \{x_0\}\bigr).
\end{equation}
При этом имеет место нормировка
\begin{equation}\label{v000o}
	\lim_{\widetilde{D}\ni x\to x_0}\frac{\widetilde{v}(x)}{-h_m\bigl(|x-x_0|\bigr)}\overset{\eqref{v0o}}{=} \frac{b}{a} \, ; \quad 
	\widetilde{v}\bigm|_{\Int S\setminus \{x_0\}}\overset{\eqref{v0o}}{\in} \Har \bigl(\Int S\setminus \{x_0\}\bigr).
\end{equation}
Ввиду \eqref{v000o} для  функции 
$\widetilde{V}\overset{\eqref{v000}}{:=}\frac{a}{b}\, \widetilde{v}	\in \sbh^+\bigl(\RR^m_{\infty}\setminus \{x_0\}\bigr) $
выполнено условие нормировки --- последнее  равенство в \eqref{VVV}. Положим 
	$\widetilde{c}:=\frac{a}{b} =  \frac1b \inf\limits_{x\in \partial S}  g_{\widetilde{D}}(x, x_0) $
как в \eqref{cz0}. Тогда из    \eqref{v000}, домноженного на  $\widetilde{c}$, получим в точности \eqref{cr:Vng}. Наконец, все перечисленные свойства из \eqref{VVV}  ---  прямые следствия построения \eqref{v000}--\eqref{v000o} и известных свойств функции Грина $g_{\widetilde{D}}(\cdot , z_0)$, участвующей в нём, начиная с  конструкции функции $v_0\:=\frac{b}{a}\,  g_{\widetilde{D}}(\cdot , x_0)$. 
\end{proof}

\subsubsection{{Продолженная функция $\widetilde V$ --- предел потенциалов Йенсена}}\label{sss242}
\begin{propos}\label{pr:Vn}
Пусть $\widetilde{V}$ --- функция из \eqref{cr:Vng}, построенная в\/
{\rm пп.~\ref{spf:241}.} Тогда 
\begin{equation}\label{Vninc+}
{V}_n:=\Bigl(\widetilde{V}-\frac1n\Bigr)^+:=\max\Bigl\{0, \widetilde{V}-\frac1n\Bigr\}\in 	\sbh^+\bigl(\RR^m_{\infty}\setminus\{x_0\}\bigr), \quad n\in \NN,
\end{equation}
--- потенциалы  Йенсена $V_n\in PJ_{x_0}(D)$,   для которых поточечно 
$\lim\limits_{n\to+\infty} V_n=\widetilde{V}$, $V_n\leq V_{n+1}$   на $D\setminus \{x_0\}$ при всех  $n\in \NN$.
Кроме того, для некоторого числа $r_0>0$ с вложением $B_*(x_0,r_0) :=B(x_0,r_0)\setminus \{x_0\}\Subset S$ и для некоторого номера  $n_0\in \NN$ имеем
\begin{equation}\label{cr:Vnh}
V_n\bigm|_{B_*(x_0,r_0)}\in \Har \bigl(B_*(x_0,r_0)\bigr)\quad \text{при всех  $n\geq n_0$};
\quad \lim_{(D\setminus \{x_0\})\ni x \to x_0}
\frac{V_n(x )}{-h_m\bigl(|x-x_0|\bigr)}\overset{\eqref{nvz}}{=} 1.
\end{equation}
\end{propos}

\begin{proof}
При любом $n\in \NN$ для каждой функции $V_n$ из \eqref{Vninc+}
ее субгармоничность и  положительность в $\RR^m_{\infty}\setminus \{x_0\}$ --- следствие из  \eqref{cr:Vng}, 
нормировка-равенство в  \eqref{cr:Vnh} --- из последнего равенства в \eqref{VVV},  гармоничность в  $B_*(x_0,r_0)$ в \eqref{cr:Vnh} ---  из первого соотношения в \eqref{VVV}, а требуемая в определении \ref{df:PoJ} потенциалов Йенсена финитность в $D$ вытекает из построения \eqref{Vninc+} и свойства $\lim\limits_{\partial D} \widetilde V=0$ из \eqref{VVV}. По определению \ref{df:PoJ} функции   ${V}_n \in PJ_{x_0}(D)$ --- потенциалы Йенсена внутри $D$ с полюсом $x_0\in D$. По построению \eqref{Vninc+} последовательность $(V_n)_{n\in \NN}$ возрастающая и стремится поточечно к  $\widetilde{V}$.
\end{proof}

\subsection{Доказательство основной теоремы}\label{ss:+} Будет использована построенная в  п.~\ref{sss242}
 возрастающая к функции $\widetilde{V}$ из \eqref{cr:Vng} последовательность потенциалов Йенсена $(V_n)_{n\in \NN}$ со свойствами \eqref{cr:Vnh}. Каждому потенциалу Йенсена $V_n$ по предложению \ref{pr:1} с отображением $\mathcal P$ из \eqref{con:P} соответствует 
мера Йенсена $\mu_n$ внутри области  $D$ в точке $x_0\in D\setminus {\dom M}$: 
\begin{equation}\label{cr:Vnh00}
\mu_n:= \mathcal P^{-1} (V_n)\overset{\eqref{eq:mu}, \eqref{cr:Vnh}}{=}\frac1{s_{m-1}} \Delta V_n\in PJ_{x_0}(D),\quad 
\supp \mu_n\overset{\eqref{cr:Vnh}}{\subset} D\setminus B(x_0,r_0),
\end{equation}
 где $r_0>0$ и $B(x_0,r_0)\subset \Int S$. Пусть $M=u_1-u_2$, где $u_1,u_2\in \sbh (D)$ соответственно с мерами Рисса  $\nu_M^+, \nu_M^-\in \mathcal M^+(D)$. По определению $\delta$-субгармонической функции ввиду  $x_0\in {\dom M}$ имеем $u_1(x_0)\neq -\infty$ и $u_2(x_0)\neq -\infty$. 
При этом,  если $u(x_0)=-\infty$ или \underline{не} выполнены условия 
	\begin{equation}\label{Mv}
				\int_{D\setminus S}v {\d}  \nu_M^-<+\infty,  \quad 
					\int_{\widetilde{D}\setminus S}g_{\widetilde{D}} (\cdot, x_0){\d}  \nu_M^- 	{<}+\infty,  	
				\end{equation}
то неравенство  \eqref{mest} тривиально. Поэтому можем считать, что  $u(x_0)\neq -\infty$ и одновременно выполнены соотношения  \eqref{Mv}. По расширенной формуле Пуассона\,--\,Йенсена \eqref{f:PJ} из предложения \ref{pr:2}, применённой к субгармоническим функциям $u, u_1,u_2$,   получаем
\begin{subequations}\label{PJ}
\begin{align}
	u(x_0) +\int_{D\setminus \{x_0\}} V_{n} {\d}  {\nu}_u &\overset{\eqref{f:PJ}, \eqref{cr:Vnh00}}{=}\int_{D\setminus B(x_0,r_0)} u {\d}  \mu_n  ,
\tag{\ref{PJ}$u$}\label{Pju}\\ 
u_1(x_0) +\int_{D\setminus \{x_0\}} V_{n} {\d}  {\nu}_M^+ &\overset{\eqref{f:PJ},\eqref{cr:Vnh00}}{=}\int_{D\setminus B(x_0,r_0)} u_1 {\d}  \mu_n  ,
\tag{\ref{PJ}$u_1$}\label{Pju1} \\
u_2(x_0) +\int_{D\setminus \{x_0\}} V_{n} {\d}  {\nu}_M^- &\overset{\eqref{f:PJ}, \eqref{cr:Vnh00}}{=}\int_{D\setminus B(x_0,r_0)} u_2 {\d}  \mu_n  .
\tag{\ref{PJ}$u_2$}\label{Pju2}
\end{align}
\end{subequations}
Из условия $u\leq M=u_1-u_2$ на $D$ для правых частей равенств \eqref{PJ} получаем
\begin{equation*}
	\int_{D\setminus B(x_0,r_0)} u {\d}  \mu_n \leq \int_{D\setminus B(x_0,r_0)} M {\d}  \mu_n=
	\int_{D\setminus B(x_0,r_0)} u_1 {\d}  \mu_n- \int_{D\setminus B(x_0,r_0)} u_2 {\d}  \mu_n.
\end{equation*}
Отсюда  по трём  равенствам \eqref{PJ}
\begin{multline}\label{uu122}
	u(x_0) +\int_{D\setminus \{x_0\}} V_{n} {\d}  {\nu}_u 	+\int_{D\setminus \{x_0\}} V_{n} {\d}  {\nu}_M^-
\\	\leq M(x_0) +\int_{D\setminus \{x_0\}} V_{n} {\d}  {\nu}_M^+
		\leq M(x_0) +\int_{D\setminus \{x_0\}} \widetilde{V} {\d}  {\nu}_M^+,
\end{multline}
поскольку  последовательность функций $V_n$, возрастает и стремится к $\widetilde{V}$ на $D$ поточечно. Для интеграла в  правой части допускается и значение $+\infty$. Если этот интеграл действительно равен $+\infty$, то ввиду конечности интегралов \eqref{dMa} для $x=x_0\in {\dom M}$ получаем
$+\infty=\int_{D\setminus \widetilde{D}} \widetilde{V} {\d}  {\nu}_M^+\overset{\eqref{cr:Vng}}{=}
	\int_{D\setminus \widetilde{D}} \widetilde{c} \, v {\d}  {\nu}_M^+$.
Отсюда, ввиду конечности первого интеграла в \eqref{Mv}, первый  интеграл в правой части   \eqref{mest} также равен $+\infty$ и доказывать нечего. Поэтому далее предполагаем, что интеграл в правой части  \eqref{uu122} конечен. Применяя теорему 
 о монотонной сходимости для интегралов  к левой части \eqref{uu122}, получаем 
\begin{equation}\label{itots}
	u(x_0) +\int_{D\setminus \{x_0\}} \widetilde{V} {\d}  {\nu}_u 	+\int_{D\setminus \{x_0\}} \widetilde{V}{\d}  {\nu}_M^-
	\leq M(x_0) +\int_{D\setminus \{x_0\}} \widetilde{V} {\d}  {\nu}_M^+.
\end{equation}
Здесь  в силу \eqref{cr:Vng}
\begin{equation*}
		\widetilde{V}= 
	\begin{cases}
	g_{\widetilde{D}}(\cdot, x_0)\quad &\text{на $S\setminus \{x_0\}$},\\
	\max\bigl\{g_{\widetilde{D}}(\cdot, x_0),\;\widetilde{c}\, v \bigr\}\quad &\text{на $\widetilde{D}\setminus S$},\\
	\widetilde{c}\, v \quad &\text{на $D\setminus \widetilde{D}$}.
	\end{cases}
	\end{equation*}
		Из этих равенств применительно к \eqref{itots}  получаем
		\begin{multline}\label{es:uM}
		u(z_0) +\int_{S\setminus \{x_0\}} g_{\widetilde{D}}(\cdot, x_0)  {\d}  {\nu}_u  
		+\int_{D\setminus S} \widetilde{c}\, v {\d}  {\nu}_u 	\\		\leq 
		M(x_0) +\int_{S\setminus \{x_0\}} g_{\widetilde{D}}(\cdot, x_0)  {\d}  {\nu}_M  
		+\int_{\widetilde{D}\setminus S} \max\bigl\{g_{\widetilde{D}}(\cdot, x_0),\;\widetilde{c}\, v \bigr\}  {\d}  {\nu}_M  		+\int_{D\setminus \widetilde{D}} \widetilde{c}\, v {\d}  {\nu}_M.
				\end{multline}
Второе положительное слагаемое-интеграл в левой части \eqref{es:uM} можно убрать.  Деление   \eqref{es:uM} на  $\widetilde{c}>0$ из  \eqref{cz0} даёт  с $C :=\dfrac{1}{\widetilde{c}}\overset{\eqref{cz0}}{=}\dfrac{b} {\inf\limits_{x\in \partial S}  g_{\widetilde{D}}(x, x_0)}  >0$ вида  \eqref{cz0C} неравенство
\begin{multline*}
C u(x_0) 	+\int_{D\setminus S}  v {\d}  {\nu}_u 	\leq 	\int_{D\setminus \widetilde{D}}  v {\d}  {\nu}_M	
\\	+		C\int_{S\setminus \{x_0\}} g_{\widetilde{D}}(\cdot, x_0)  {\d}  {\nu}_M  
		+\int_{\widetilde{D}\setminus S} \max\bigl\{Cg_{\widetilde{D}}(\cdot, x_0),\; v \bigr\}  {\d}  {\nu}_M^+  +CM(x_0).
\end{multline*}	
Отсюда в силу очевидного для положительных функций  неравенства
$\max\bigl\{Cg_{\widetilde{D}}(\cdot, x_0),\; v \bigr\}\leq Cg_{\widetilde{D}}(\cdot, x_0)
	+ v$  на $D\setminus S$
 получаем 
\begin{multline*} 
C u(x_0) 	+\int_{D\setminus S}  v {\d}  {\nu}_u 	\leq 	\int_{D\setminus S}  v {\d}  {\nu}_M	+\int_{\widetilde{D}\setminus S}
v {\d}  {\nu}_M^-
\\	+		C\int_{\widetilde{D}\setminus \{x_0\}} g_{\widetilde{D}}(\cdot, z_0)  {\d}  {\nu}_M  
		+C\int_{\widetilde{D}\setminus S} g_{\widetilde{D}}(\cdot, z_0)  {\d}  {\nu}_M^-  +CM(x_0),
		 \end{multline*}	
что и доказывает требуемое \eqref{mest}. 	Наконец, при $\widetilde{D}\Subset D$  как второй интеграл в правой части \eqref{mest} ввиду ограниченности $v$ на $\widetilde{D}\setminus S$, так и интегралы, участвующие в определении величины $\overline{C}_M$,  для $x_0\in \dom M$ ввиду \eqref{dMa} с $x=x_0$ конечны. Основная теорема доказана.
	
	\subsection{Субгармоническая мажоранта $M$} Всюду ниже $D\subset \CC_{\infty}^n$
	и $z_0=x_0\in \Int S \Subset D$. 
 \begin{proof}[ теоремы \ref{th:1}] Очевидно, можно считать, что $f(z_0)\neq 0$ и интеграл в правой части
\eqref{iq:ufu} конечен. По основной теореме при $u:=\log |f|\in \sbh_*(D)$ для частного случая  $ M \in \sbh_* (D)$  основной теоремы, очевидно, 
$\nu_M^-=0$  и в правой  части  \eqref{mest}, как и выражении для $\overline C_M$, интегралы по мере $\nu_M^-$  равны нулю. 
Всегда можно подобрать  регулярную  область $\widetilde{D}\Subset D$, описанную в   \cite[1.4.3, определение]{HK}, для которой $S\Subset \widetilde{D}$. Таким образом, для любой тестовой  функции $v\overset{\eqref{sbh+}}{\in} \sbh_0^+(D\setminus S;\leq b)$ из неравенства \eqref{mest} получаем
\begin{equation}\label{ittCMg+}
\int_{D\setminus S}  v \d {\nu}_u 	\leq 	\int_{D\setminus S}  v \d  {\nu}_M	 - C u(z_0) 	
+		C\int_{\widetilde{D}} g_{\widetilde{D}}(\cdot, z_0)  {\d}  {\nu}_M  		 +CM^+(z_0).
		 \end{equation}	
Выбор области $\widetilde{D}$ полностью обусловлен лишь взаимным расположением  $S\Subset D$ и области  $D$, т.\,е. в выборе постоянной $C$ из \eqref{cz0C} влияние области $\widetilde{D}$ можно заменить на зависимость от  $S$ и $D$. Из этих же соображений  постоянная  
	$\overline{C}_M\overset{\eqref{mest}}{:=} \int_{\widetilde{D}} g_{\widetilde{D}}(\cdot, z_0)  {\d}  {\nu}_M   +M^+(z_0)$
зависит только $z_0$, $S$, $D$ и $M$, как и требуется. 
 По формуле Пуанкаре\,--\,Лелона из пп.~\ref{zf} левую часть в \eqref{ittCMg+} можно заменить на левую часть неравенства \eqref{iq:ufu}. Теорема \ref{th:1} доказана.
\end{proof}
 \begin{proof}[ теоремы \ref{th:2}] Исходя из теоремы \ref{th:1}, достаточна 
 \begin{lemma}\label{pr:tmrw}
Для   полунепрерывной сверху функции $w$ из \eqref{c:w} при любом из двух условий\/ 
{\rm \eqref{key}(i)} или\/ {\rm \eqref{key}(ii)} полунепрерывная  сверху регуляризация $\gm^*w$ 
её  наибольшой миноранты $\gm w$ относительно $\sbh_0(D\setminus S)$ --- тестовая функция 
из  $\sbh_0^+(D\setminus S;<+\infty)$, удовлетворяющая ограничению $\gm^*w\leq w$ на $D\setminus S$.
 \end{lemma}	
\begin{proof}[леммы \ref{pr:tmrw}] Даже без условий \eqref{key}  по построению 
$\gm^* w\in \sbh^+(D\setminus S)$ и $\gm^* w\leq w$ на $D\setminus S$ ввиду полунепрерывности 
сверху функции $w$. Тогда при условии \eqref{key}(i)
имеем $0\leq \lim\limits_{\partial D} \gm^* w \leq \lim\limits_{\partial D} w=0$ и $\gm^* w 
\overset{\eqref{sbh+S}}{\in} \sbh_0^+(D\setminus S;<+\infty)$, что и требуется.

В случае \eqref{key}(ii)  пусть $z_0\in \Int  S$ и $B(z_0,r_0)\subset  S$ для некоторого $r_0>0$. 
Используем продолженную на $\complement D$ 
функцию Грина $g_D(\cdot, z_0)$  из  \eqref{gD0r}, где в роли  $ \widetilde D$ --- область  $D$. 
Из регулярности области $D$ следует $g_D(\cdot ,z_0)=0$ на $\partial D$, а ввиду 
$B(z_0,r_0)\subset  S\Subset D$ найдётся такое $C\in \RR^+$, что $Cg_D(\cdot,z_0 )
\overset{\eqref{c:w}}{\geq}	b$ на $\partial S$. При этом $Cg_D(\cdot, z_0)\bigm|_{D\setminus \{z_0\}}\in \Har 	\bigl(D\setminus \{z_0\}\bigr)$, $Cg_D(\cdot, z_0)\in \sbh	\bigl(\CC_{\infty}^n\setminus \{z_0\}\bigr)$. В частности,  продолженная функция Грина 
	$g_D(\cdot, z_0)$ полунепрерывна сверху на $\CC_{\infty}^n\setminus \{z_0\}$.
	Для любой функции $v\in \sbh_0( D\setminus S)$, удовлетворяющей ограничению $v\leq w$ 
	на $D\setminus S$, по принципу максимума имеем $Cg_D(\cdot , x_0){\geq} v$ на $D\setminus {S}$ ввиду гармоничности 	левой части в $D\setminus \{z_0\}$. Отсюда  и из  полунепрерывности сверху функции $Cg_D(\cdot, z_0)$ на 	$\CC_{\infty}^n\setminus \{z_0\}$ следует $\gm^* w\leq Cg_D(\cdot, z_0)$ на $ D\setminus {S}$, что даёт $\lim\limits_{\partial D} \gm^* w=0$ и $\gm^* w\overset{\eqref{sbh+b}}{\in} \sbh_0^+(D\setminus S; \leq b)$. Лемма \ref{pr:tmrw}   доказана.
\end{proof}
\end{proof}

\begin{proof}[теоремы \ref{cor:1}]	 	
	Пусть $v\overset{\eqref{sbh+par}}{\in}\sbh_0^+(\partial D)$ и $v\leq w$ на некотором борелевском  множестве $B$ из $\Conc \nu_M$. Тогда по определению \ref{testv}  тестовых функций можно подобрать  точку $z_0\in D$, для которой $M(z_0)\neq -\infty$ и $f(z_0)\neq 0$, а также замкнутое подмножество $S\Subset D$ с $\Int S\ni z_0$ так, что 	функция $v$ определена как субгармоническая на $D\setminus \Int S$ и   $v\overset{\eqref{sbh+b}}{\in} \sbh_0^+(D\setminus S;\leq b)$ для $b=\sup\limits_{\partial S} v<+\infty$, а $B_S:=B\cap (D\setminus S) \in \mathcal B (D\setminus S)$
	и $\int_{B_S}w \d \nu_M<+\infty$. По основной  теореме с $u:=\log |f|$  неравенство  \eqref{mest} с потерей, как и при доказательстве теоремы \ref{th:1},  части информации  можно записать для некоторой постоянной $C_0\in \RR^+$ в виде 
$\int_{D\setminus {S}}  v \d {\nu}_u \overset{\eqref{ittCMg+}}{\leq}	\int_{D\setminus S}  v \d {\nu}_M	+C_0= 	\int_{B_S}  v \d {\nu}_M	+C_0\leq 	\int_{B_S}  w \d {\nu}_M	+C_0<+\infty$.  По формуле Пуанкаре\,--\,Лелона из пп.~\ref{zf} в обозначении \eqref{in_p} левую часть
здесь можно заменить на интеграл  $ \int^{\partial D} v \, {\rm Z} \d \sigma_{2n-2}\overset{\eqref{in_p}}{<} +\infty$
	для любого поддивизора ${\rm Z} \leq  \Zero_f$. Теорема \ref{cor:1} доказана.
\end{proof}
\paragraph{Заключительные замечания} {\bf 1.} Разнообразные способы построения  конкретных   тестовых функций из определения \ref{testv} и их подклассов (радиальных для круга, плоскости, кольца, выражающихся через функции Грина, функцию расстояния до подмножества на $\partial D$, гиперболический, или конформный, радиус области и т.\,д.)  приведены в \cite{KhT16} для одной комплексной переменной. В \cite{LKhT} даны начальные конструкции тестовых функций (радиальных для шара и выражающихся через функцию Грина) для многих переменных. 

{\bf 2.} <<Плюрисубгармонический>>  переход на псевдовыпуклые области в $\CC^n$ для плюрисубгармонической мажоранты  $M$  с  существенным расширением класса тестовых функций при $n>1$  несколько   сложнее, но  всё же  вполне  перспективен на основе комплексной теории потенциала с привлечением техники аналитических дисков \cite{Klimek}, \cite{Khab01},  \cite{Khsur} и теории голоморфных потоков и дисковых огибающих \cite{Po99}.

\end{fulltext}

\end{document}